\documentclass[12pt]{article}
\usepackage[final]{epsfig}
\usepackage{graphics}
\usepackage{amsmath}
\usepackage{amsfonts}
\usepackage{latexsym}
\usepackage{amssymb}
\usepackage{graphicx}
\usepackage{epstopdf}
 
\newtheorem{lemma}{Lemma}[section]
\newtheorem{proposition}[lemma]{Proposition}
\newtheorem{remark}[lemma]{Remark}

\newtheorem{theorem}[lemma]{Theorem}

\newtheorem{corollary}[lemma]{Corollary}

\newtheorem{conjecture}[lemma]{Conjecture}

\begin{document}
\newcommand{\eps}{{\varepsilon}}
\newcommand{\proofend}{$\Box$\bigskip}
\newcommand{\C}{{\mathbf C}}
\newcommand{\Q}{{\mathbf Q}}
\newcommand{\R}{{\mathbf R}}
\newcommand{\Z}{{\mathbf Z}}
\newcommand{\RP}{{\mathbf {RP}}}
\newcommand{\CP}{{\mathbf {CP}}}
\newcommand{\A}{{\rm Area}}
\newcommand{\Le}{{\rm Length}}


\newcommand{\marginnote}[1]
{
}

\newcounter{ml}
\newcommand{\bk}[1]
{\stepcounter{ml}$^{\bf ML\thebk}$%
\footnotetext{\hspace{-3.7mm}$^{\blacksquare\!\blacksquare}$
{\bf ML\thebk:~}#1}}

\newcounter{st}
\newcommand{\st}[1]
{\stepcounter{st}$^{\bf ST\thest}$%
\footnotetext{\hspace{-3.7mm}$^{\blacksquare\!\blacksquare}$
{\bf ST\thest:~}#1}}


\title {On bicycle tire tracks geometry,  hatchet planimeter, Menzin's conjecture and oscillation of unicycle tracks}
\author{Mark Levi\thanks{
Department of Mathematics,
Pennsylvania State University, University Park, PA 16802, USA;
e-mail: \tt{levi@math.psu.edu. }
}
\, and Serge Tabachnikov\thanks{
Department of Mathematics,
Pennsylvania State University, University Park, PA 16802, USA;
e-mail: \tt{tabachni@math.psu.edu}
}
\\
}

\date{\today}
\maketitle
\begin{abstract} 
The model of a bicycle is a unit segment $AB$ that can move in the plane so that it remains tangent to the trajectory of point $A$ (the rear wheel is fixed on the bicycle frame); the same model describes the hatchet planimeter. The trajectory of the front wheel  and the initial position of the bicycle uniquely determine  its motion and its terminal position; the monodromy map sending the initial position to the terminal one arises. According to R. Foote's theorem, this mapping of a circle to a circle is a Moebius transformation. We extend this result to multi-dimensional setting. Moebius transformations belong to one of the three types: elliptic, parabolic and hyperbolic. We prove a 100 years old Menzin's conjecture: if the front wheel track is an oval with area at least $\pi$ then the respective monodromy is hyperbolic. We also study bicycle motions introduced by D. Finn in which the rear wheel follows the track  of the front wheel. Such a ''unicycle" track becomes more and more oscillatory in forward direction. We prove that it cannot be infinitely extended backward and relate the problem to the geometry of  the space of forward semi-infinite equilateral linkages.
\end{abstract}

\section{Introduction} \label{intro}

The geometry of bicycle tracks is a rich and fascinating subject. Here is a sampler of questions:
\begin{enumerate}
\item Given the tracks of the rear and front wheel, can you tell which way the bicycle went?

\item The track of the front wheel is a smooth simple closed curve. Can one ride the bicycle so that the  rear wheel track also closes up?

\item Can one ride a bicycle in such a way that the tracks of the rear and front wheels coincide\footnote{Other than along a straight line.}?
\end{enumerate}

Our model of a bicycle is an oriented segment, say, $AB$, of length $\ell$ that can move in the plane in such a way that the trajectory of point $A$  always remains tangent to the segment.  Point $A$ represent the rear wheel, point $B$ the front wheel; the rear wheel is fixed on the bicycle frame whereas the front wheel can turn, and this explains the law of motion. (Most often we set $\ell=1$ -- this can be always assumed by making a dilation -- but sometimes we shall consider $\ell$ as a parameter and allow it to take very small or very large values.) Thus  the end point of the oriented tangent segment to the trajectory of the rear wheel traces the trajectory of the front wheel, see \cite{Fi,KVW}.

The same mathematical model describes another mechanical device, the Prytz or hatchet planimeter, see \cite{Ba,Cr,Fo}.  Various kinds of planimeters were popular objects of study in the late 19th and early 20th century.

The first of the above questions has the following answer: generically, one  can determine the direction, but in some special cases one cannot: for example, for concentric circles of radii $r$ and $R$ satisfying $r^2+\ell^2=R^2$. Surprisingly, the problem of describing such ``ambiguous" pairs of closed tracks is equivalent to Ulam's problem of describing (2-dimensional) bodies that float in equilibrium in all positions. See \cite{Ta,We,We1,We2} for a variety of results and  references.

The content of the present paper has to do with the other two questions. In Section \ref{prel} we place the problem into the framework of contact geometry. We allow the trajectory of the rear wheel to be a wave front, that is, to have cusp singularities, but we show that the trajectory of the front wheel remains smooth. We  deduce a useful differential equation relating the motions of the rear and the front wheels.

\begin{figure}[hbtp]
\centering
\includegraphics[width=3in]{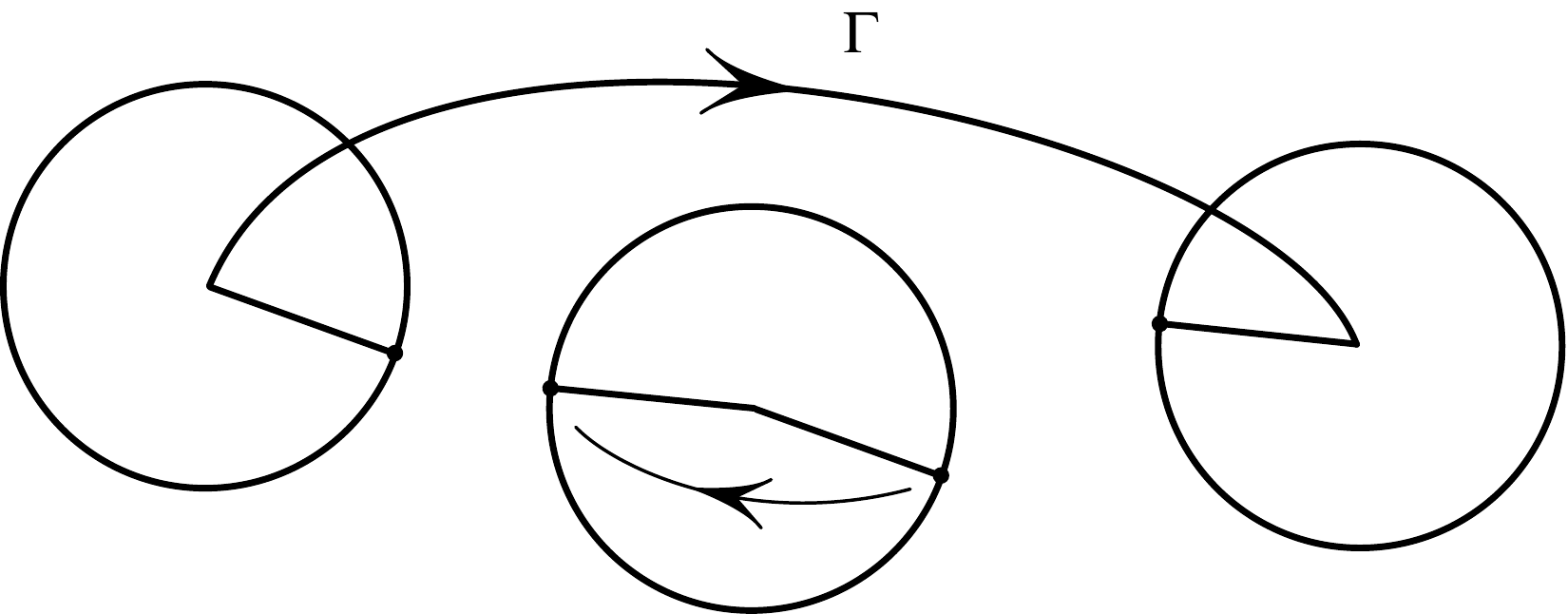}
    	\caption{The circle mapping generated by the curve $\Gamma$. According to Foote's theorem, this mapping is a  M\"obius transformation.}
    	\label{fig:moebius}
\end{figure}

Fixing a path $\Gamma$ of the front wheel gives rise to a circle map: the initial direction of the segment, characterized by a point on the circle, determines its final direction, 
see figure~\ref{fig:moebius}. We will refer to this map of the circle to itself (the two circles are identified by parallel translation) as the monodromy map.\footnote{T. Tokieda suggested the term ``opisthodromy" (literary, rear track).}   It is a beautiful theorem of R. Foote
  \cite{Fo} (see also \cite{L-W}) that, for every trajectory  of the front wheel, the 
  monodromy map is a M\"obius transformation.  
  In Section \ref{monod} we reprove this theorem and extend
   it to bicycle motion in Euclidean space of any dimension.
  
A (non-trivial) M\"obius transformation is of one of the three types: elliptic, parabolic or hyperbolic. The former have no fixed points and the latter have exactly two, one attracting and one repelling (parabolic ones have a single, neutral fixed point). Suppose the trajectory of the front wheel is a closed curve. Then, up to conjugation, the respective monodromy, and therefore its type, does not depend on the initial point. In Section \ref{monod} we give a necessary and sufficient condition for the monodromy to be parabolic: the trajectory of the rear wheel is a closed wave front with the total algebraic arc length equal to zero (the sign of the arc length changes when passing a cusp).

\begin{figure}[hbtp]
\centering
\includegraphics[width=5.5in]{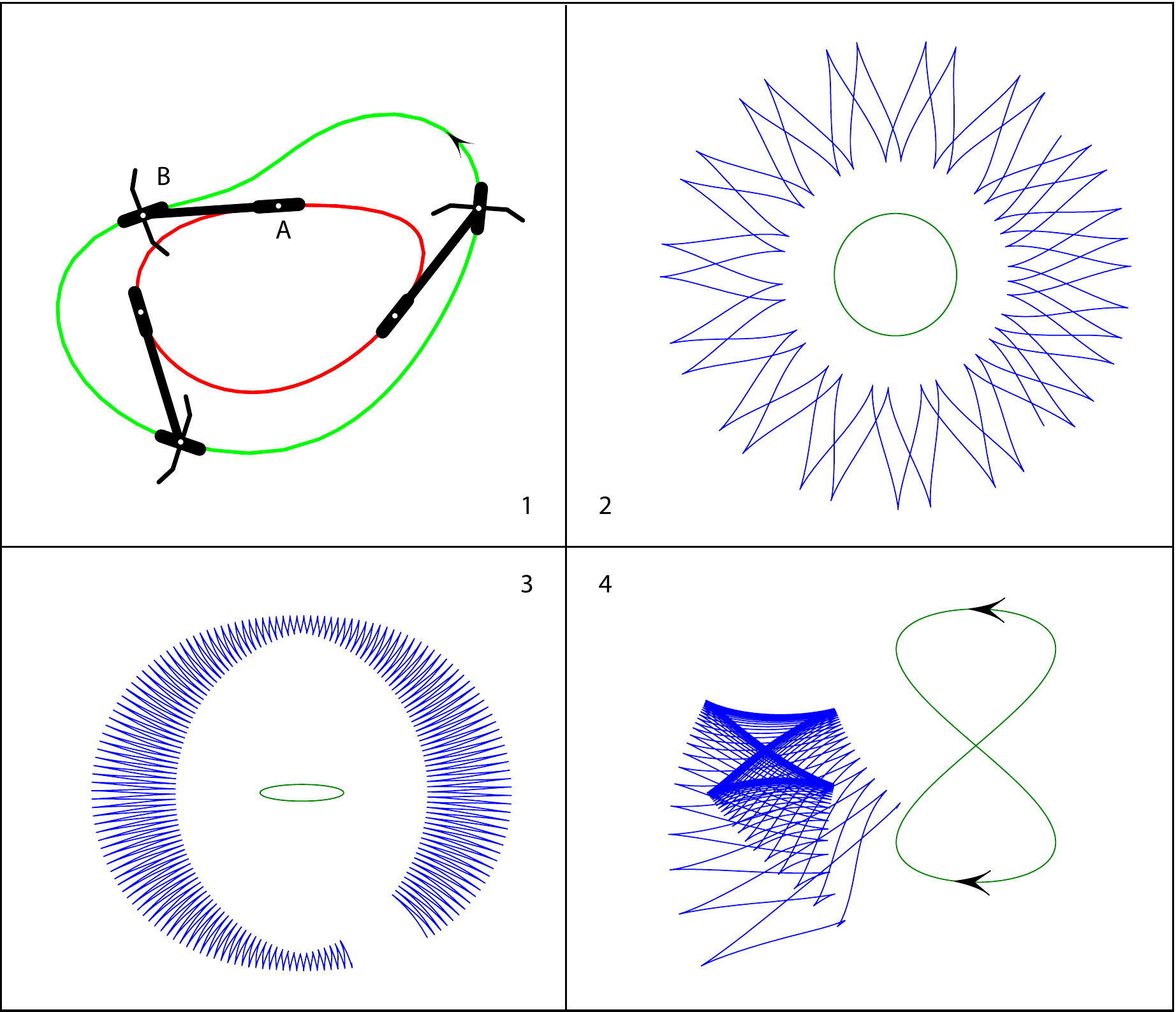}
    	\caption{Examples 1 and 4 are hyperbolic; 2 and 3 are elliptic.
	The areas bounded by the two curves in 1 differ by $ \pi \ell ^2 $. }
    	\label{fig:examples}
\end{figure}

Still assuming that the trajectory of the front wheel is  closed, a fixed point of the monodromy map corresponds to a closed trajectory of the rear wheel. Thus, in the hyperbolic case, for a given closed trajectory of the front wheel, 
there are exactly two bicycle motions such that the trajectory of the rear wheel is closed; each of these motions is hyperbolically attracting for one of the choices of the direction of motion; examples are shown in figure~\ref{fig:examples} (1) and (4). In contrast, in the elliptic case, no trajectory of the rear wheel closes after one cycle. It is worth mentioning that, for some trajectories of the front wheel, the monodromy is the identity: for every bicycle motion the trajectory of the rear wheel closes up.

A hundred years old conjecture by Menzin \cite{Me} states, in our terms, that if the trajectory of the front wheel is a closed convex curve bounding area greater than $\pi \ell^2$, then the respective monodromy is of the hyperbolic type. In Section \ref{Menzin} we prove this conjecture. The main tool is the classical Wirtinger inequality. Earlier Foote \cite{Fo} proved Menzin's conjecture for parallelograms.

\begin{figure}[hbtp]
\centering
\includegraphics[width=5.5in]{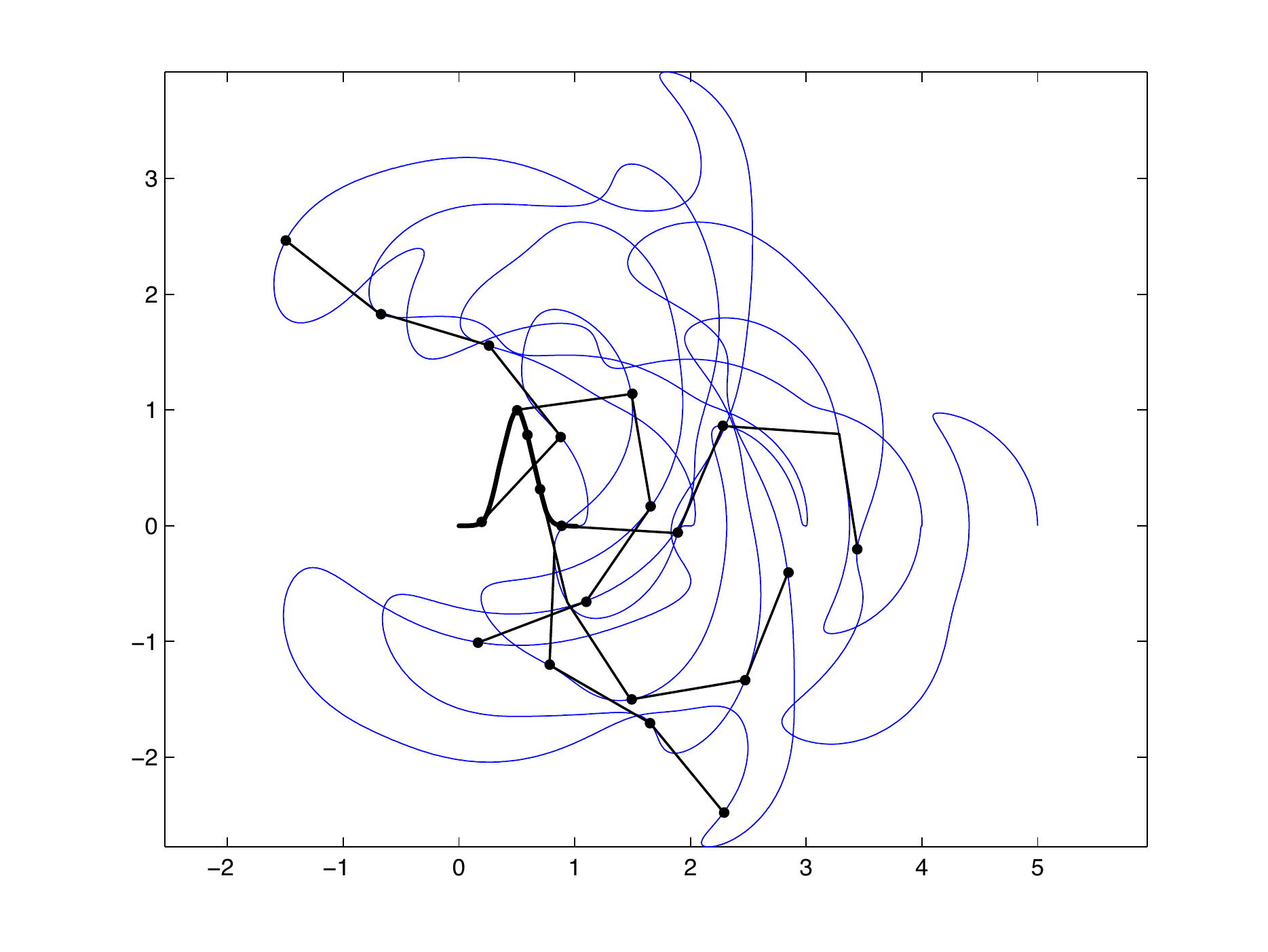}
    \caption{The rear wheel follows the track of the front wheel. The ``seed" curve is shown in heavier stroke. Several consecutive position of the associated moving linkage are shown. The shape of the linkage is very sensitive to the position of the starting point   on the seed curve.   }
    \label{fig:unicycle}
\end{figure}

Section \ref{osc} concerns Finn's construction of bicycle motion leaving a single track \cite{Fi}. Consider a ``seed" curve, tangent to the $x$-axis at  points $0$ and $1$ with all derivatives and oriented to the right (the ``fat" curve in figure~\ref{fig:unicycle}). This curve is the initial trajectory of the rear wheel; drawing the tangent segments of length 1 to it yields the next curve which is tangent to the $x$-axis at  points $1$ and $2$ with all derivatives. Iterating this process, one obtains  bicycle motion leaving a unicycle track, i.e., a curve which both wheels follow. 

Numerical study shows that, unless the seed curve is horizontal, the resulting unicycle track becomes more and more oscillating, figure ~\ref{fig:unicycle}. We prove that the number of intersections with the $x$-axis and the number of extrema of the height function increase at least by one with every iteration of this construction.  As a consequence, the seed curve with finitely many intersections with the $x$-axis (or a finite number of extrema) has at most finitely many preimages under Finn's construction. This means that the corresponding unicycle track cannot extend back indefinitely. We also make a number of conjectures on the Finn construction strongly supported by numerical evidence.

A unicycle track can be viewed as an integral curve of a direction field in a certain infinite dimensional space. Specifically, we consider  the configuration space of equilateral forward infinite linkages in the plane. We  constrain the velocity of the  
$i$th vertex to the direction of the $i$th link (heuristically, the $i$th link is the position of the bike on the $i-1$st step of Finn's construction).  This constraint defines a field of directions. Now, a forward bicycle motion generating a single track corresponds to a particular integral curve of this field of directions. This field does not satisfy the uniqueness property: through every points there pass infinitely many smooth integral curves. We also generalize Finn's construction for an arbitrary initial equilateral forward infinite linkage in which the adjacent links are not perpendicular (the Finn construction corresponds to a linkage aligned along a line).

\bigskip

{\bf Acknowledgments}. It is a pleasure to thank M. Kapovich, R. Montgomery, A. Novikov, R. Schwartz, S. Wagon and V. Zharnitsky for their interetst and help. The first author was supported by an NSF grant   DMS-0605878 and the  second one  by an NSF grant DMS-0555803.

\begin{figure}[hbtp]
\centering
\includegraphics[width=5in]{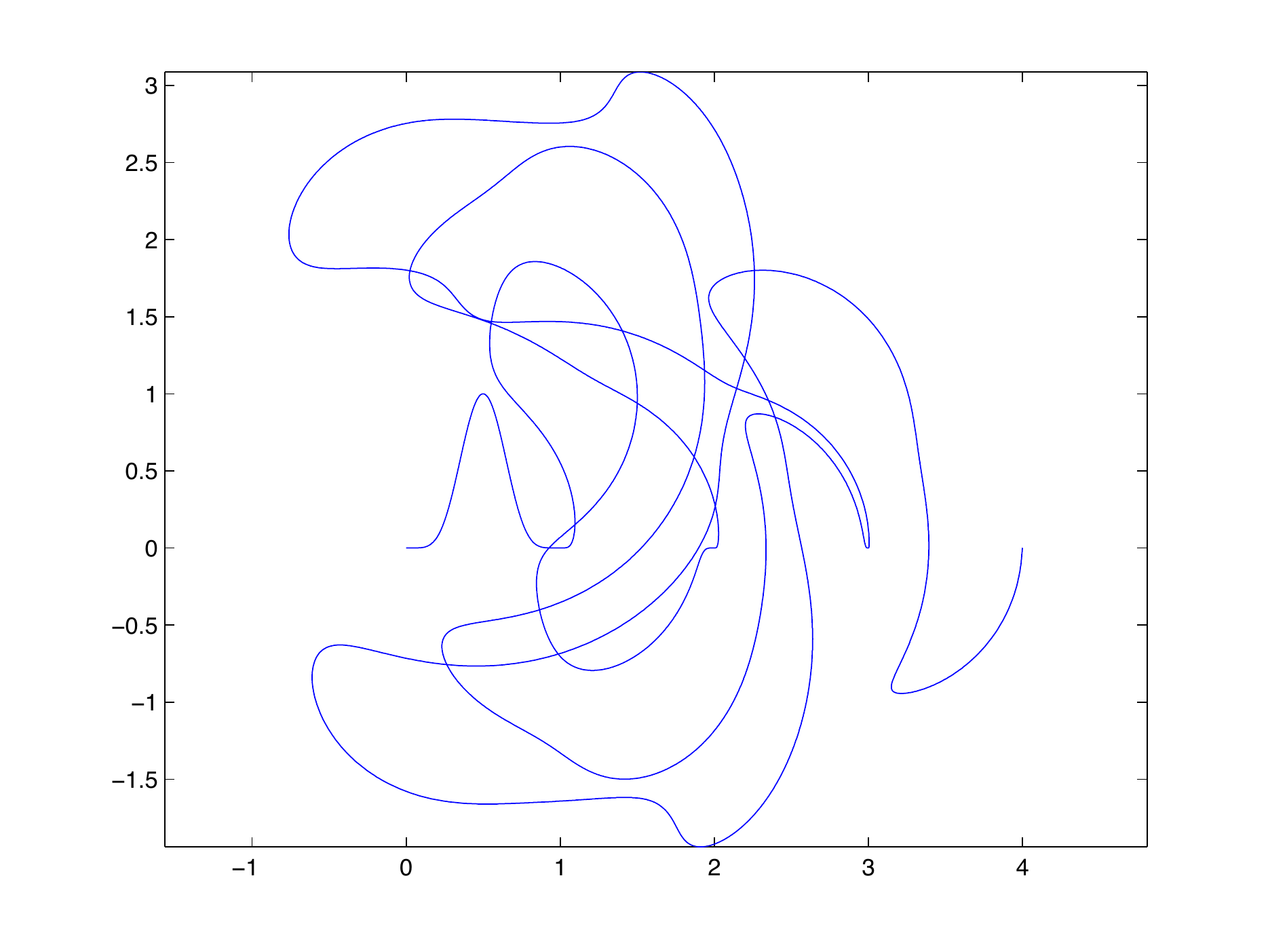}
    	\caption{The figure shows first four iterates of the initial seed curve
  $ y = 4^6 x^6(1-x)^6 $. Since this curve has only a finite order of contact with the $x$--axis, 
  only finitely many iterations are defined. }
    	\label{wiggles}
\end{figure}

\section{Preliminaries: contact geometric point of view} \label{prel}

\begin{figure}[hbtp]
\centering
\includegraphics[width=1.5in]{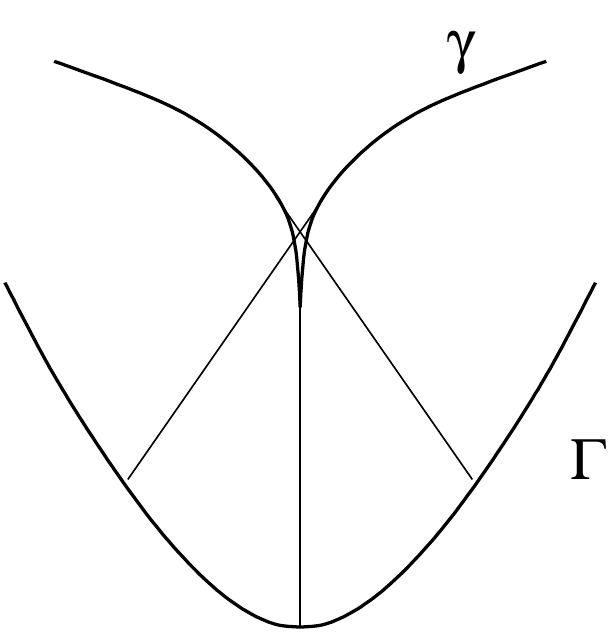}
\caption{Cusp of the curve $\gamma$}
\label{cusp}
\end{figure}

We use the notation from Section \ref{intro}.
Denote the trajectory of the rear wheel $A$ by $\gamma$ and that of the front wheel $B$ by $\Gamma$. We allow $\gamma$ to have cusp singularities as in figure \ref{cusp}. A proper prospective is provided by contact geometry, see \cite{A-G} or \cite{Ge}.

The position of the segment $ AB $ is determined by its foot point $ A(x,y) $ and by the angle $\theta$ between the $x$--axis and the segment. The infinitesimal motions in the configuration space 
$ \{(x,y, \theta) \} $ are restricted by the non--skidding condition:  
 $ ( \dot x, \dot y  ) \parallel ( \cos \theta , \sin \theta )$. This condition defines a field of 
 tangent $2$-planes in the configuration space.  
This field of planes  is non--integrable and is defined by the contact  
$ 1 $-form   $\lambda=\sin\theta\ dx - \cos\theta\ dy$. 

A smooth curve in a contact manifold is called Legendrian if its tangent line at every point 
lies in the contact plane. Let $p:M\to \R^2$ be the projection taking a contact element to 
its foot point. The image of a Legendrian curve is called a wave front; 
generically it is a piece-wise smooth curve with semi-cubical cusp singularities.   The 
singularities occur at the points where the Legendrian curve is tangent to the fibers of the 
projection $p$. A wave front has a well defined tangent line at every point and can be 
uniquely lifted to a Legendrian curve in the space of contact elements.

The trajectory of the rear wheel uniquely determines the trajectory of the front wheel. Denote by $T$ the correspondence $\gamma \mapsto \Gamma$ which assigns to point  $x\in \gamma$ the end point of the unit  tangent segment to $\gamma$ at $x$. We assume that a continuous choice is made between the two orientations of the unit tangent segments at a point. This amounts to choosing a coorientation of $\gamma$: the frame formed by the coorienting vector and the chosen tangent vector is positive. 
When the bicycle segment is not unit and has length $\ell$, we denote by $T_{\ell}$ the respective transformation and by $\Gamma_{\ell}$ its image.

The following two lemmas address the smoothness issue.

\begin{lemma} \label{smooth}
If $\gamma$ is a $ C^k $, $ k \geq 1 $  curve, then   $\Gamma_{\ell}$ is $ C^{k-1} $ curve for all  $\ell>0$.
\end{lemma}

\paragraph{\bf Proof.} Let $\gamma$ be parametrized by its arc length $s$. By the definition, $ \Gamma(s) = \gamma(s) + \gamma ^\prime (s) $, and it remains only to make sure  that $ \Gamma ^\prime  = \gamma ^\prime + \gamma ^{\prime\prime} \not= 0 $. 
But the last two vectors are orthogonal and the first has unit length. 
\proofend

\begin{lemma} \label{smoother}
Even if $\gamma$ has cusps, the curve $\Gamma_{\ell}$ is  smooth for all  $\ell>0$.
\end{lemma}

\paragraph{\bf Proof.}
Let $p_1:M\to\R^2$ take segment $AB$ to point $B$. The correspondence $T_{\ell}$ is the composition of the Legendrian lifting of a wave front $\gamma$ and the projection $p_1$. We claim that the fibers of $p_1$ are everywhere transverse to the contact distribution on $M$. This would imply the statement of the lemma since the fibers of the projection are transverse to the Legendrian curve $p^{-1}(\gamma)$.

In terms of the coordinates in $M$, one has: $p_1(x,y,\theta)=(x+\ell\cos\theta, y+\ell\sin\theta)$. The vector field $v=\partial_{\theta}+\ell\sin\theta\ \partial_x-\ell\cos\theta\ \partial_y$ is tangent to the fibers of $p_1$. One has: $\lambda(v)=\ell$, therefore $v$ is everywhere transverse to the contact planes, and we are done.
\proofend

Let $\gamma$ be an oriented and cooriented closed wave front. The Maslov index $\mu(\gamma)$ is the algebraic number of cusps of $\gamma$; a cusp is positive if one traverses it along the coorientation and negative otherwise. 

Let $\gamma$ be an oriented and cooriented closed wave front. Denote by $\rho(\gamma)$ the rotation number, that is, the total (algebraic) number of turns  made by its tangent direction. Let $\Gamma=T(\gamma)$.

\begin{lemma} \label{rotation}
One has: $\rho(\Gamma)=\rho(\gamma) + \frac{1}{2}\mu(\gamma)$.
\end{lemma}

\paragraph{\bf Proof.} Consider the 1-parameter family of curves $\Gamma_{\ell}$. By Lemma \ref{smoother}, this is a continuous family of smooth curves, hence the rotation number is the same for all $\ell$. Consider the case of very small $\ell$.

\begin{figure}[hbtp]
\centering
\includegraphics[width=4in]{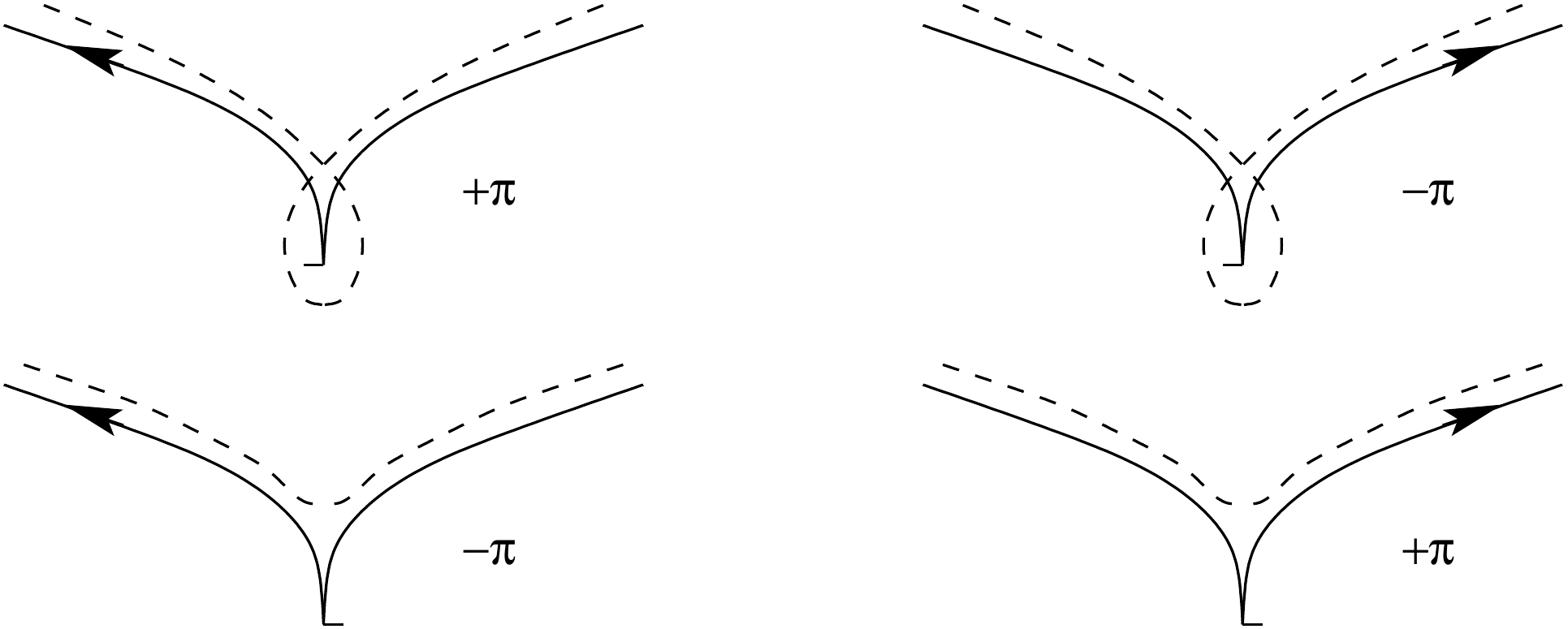}
\caption{Cusps of the curve $\gamma$ and their smoothings $\Gamma_{\ell}$}
\label{ind}
\end{figure}

Along smooth arcs of $\gamma$, the curve $\Gamma_{\ell}$ is $C^1$-close to $\gamma$. At the cusps, smoothing occurs, and the rotation of $\Gamma_{\ell}$ differs from that of $\gamma$ by $\pm \pi$. There are four cases, depending on the orientation and coorientation, depicted in figure \ref{ind}. When one traverses a cusp along the coorientation, the total rotation of  $\Gamma_{\ell}$ gains $\pi$, and when a cusp is traversed against the coorientation,  the total rotation of  $\Gamma_{\ell}$ looses $\pi$. This implies the result.
\proofend

We introduce the following notations. Let $x$ be arclength parameter along the curve $\Gamma$. The position of the segment $AB$ with $B=\Gamma(x)$ is determined by the angle made by the  tangent  vector $\Gamma'(t)$ and the vector $BA$. Let this angle be $\pi-\alpha(x)$. The function $\alpha(x)$ uniquely determines the curve $\gamma$, the locus of points $A$. Let $\kappa(x)$ be the  curvature of $\Gamma(x)$. Denote by $t$ the arclength parameter on $\gamma$ and by $k$ the curvature of $\gamma$. Note that, at cusps, $k=\infty$.

The next result is borrowed from \cite{Ta}, see also \cite{Fi}.

\begin{figure}[hbtp]
\centering
\includegraphics[width=3in]{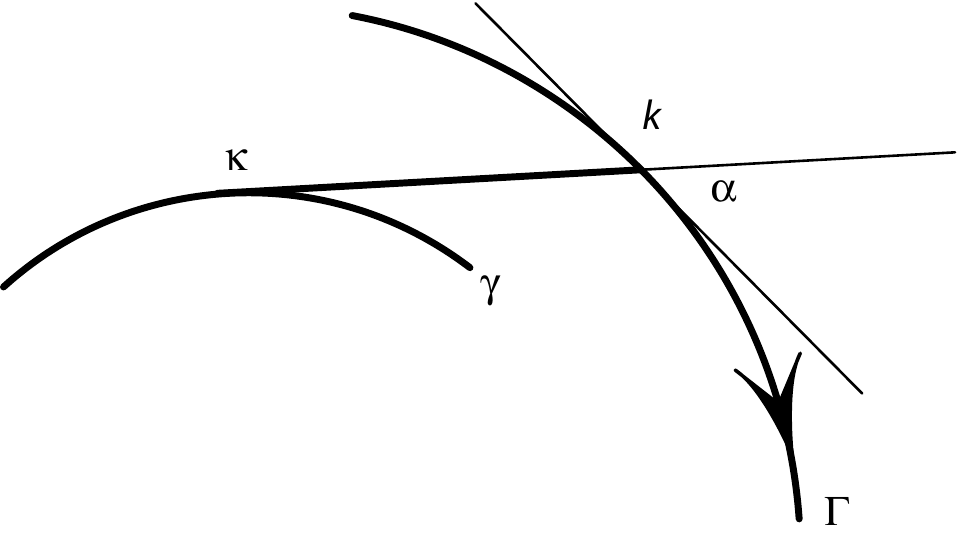}
    	\caption{Notations for Lemma \ref{local}.}
    	\label{not}
\end{figure}

\begin{lemma} \label{local}
The condition $T_{\ell}(\gamma)=\Gamma$ is equivalent to the differential equation on the function $\alpha(x)$:
\begin{equation} \label{inverse}
\frac{d\alpha(x)}{dx} + \frac{\sin \alpha(t)}{\ell}=\kappa(x).
\end{equation}
One has:
$$
\left|\frac{dt}{dx}\right|=|\cos\alpha|,\ \ k=\frac{\tan\alpha}{\ell}.
$$
\end{lemma} 

In particular the cusps of $\gamma$ correspond to the instances of $\alpha=\pm \pi/2$.

\paragraph{\bf Proof.} Let $J$ denote the rotation of the plane by angle  $\pi/2$. Then the end point of the segment of length $\ell$ making angle $\pi-\alpha(x)$ with $\Gamma'(x)$ is
\begin{equation} \label{tang}
\gamma(x)=\Gamma(x)-\ell \Gamma'(x)\cos\alpha(x)+\ell J(\Gamma'(x))\sin\alpha(x).
\end{equation}
For $T_{\ell}(\gamma)=\Gamma$ to hold, the tangent direction $\gamma'(x)$ should be collinear with the respective  segment, that is, be parallel to the vector 
$$v(x):=-\Gamma'(x)\cos\alpha(x)+J(\Gamma'(x))\sin\alpha(x).$$ 
Differentiate (\ref{tang}), taking into account that $\Gamma''(x)=\kappa(x) J(\Gamma'(x))$, and equate the cross-product with $v(x)$ to zero to obtain (\ref{inverse}).

It is straightforward to calculate that $|d\gamma/dx|=|\cos\alpha|$, hence $|dt/dx|=|\cos\alpha|$. The computation of the curvature $k$ is also straightforward.
\proofend

It is convenient to adopt the following convention: the sign of the length element $dt$ on $\gamma$ changes at each cusp. This is consistent with Lemma \ref{local} since cusps correspond to $\alpha=\pi/2$, that is, to sign changes of $\cos\alpha$. With this convention, we have $dt= \cos\alpha(x)\ dx$. In particular, the signed perimeter length of $\gamma$ is $\int_{\Gamma} \cos\alpha(x)\ dx$. 

\section{Bicycle monodromy map} \label{monod}

If $\gamma(t)$ is arc length parameterized trajectory of the rear bicycle wheel then the trajectory of the front wheel is $\Gamma(t)=\gamma(t)\pm\gamma'(t)$ (the sign depends on the coorientation of $\gamma$ and changes at its cusps). We extend this definition to bicycle rides in multi-dimensional space $\R^n$. 

On the other hand, if $\Gamma$ is given then one can recover $\gamma$, once the initial position of the bicycle is chosen. The set of all possible positions of the bicycle with a fixed position of the front wheel is a unit sphere $S^{n-1}$.  Thus there arises the time-$x$ monodromy map $M_x$ which assigns the time-$x$ position of the bicycle with a prescribed front wheel trajectory to its initial position: $M_x:S^{n-1}\to S^{n-1}$. 

 Consider the hyperbolic space $H^n$ realized as the pseudo-sphere $x_1^2+\dots + x_n^2-x_0^2=-1$ in the pseudo-Euclidean space $\R^{n,1}$ with the  metric $dx_0^2-dx_1^2-\dots -dx_n^2$. The M\"obius group $O(n,1)$ consists of linear transformations preserving the metric and acts on $H^n$ by isometries. This action extends to the null cone $x_1^2+\dots + x_n^2=x_0^2$ and to its spherization $S^{n-1}$, the sphere at infinity of the hyperbolic space. In particular, we obtain an action of the Lie algebra $o(n,1)$ on $S^{n-1}$. 
 
The following result is a multi-dimensional generalization of  Foote's theorem \cite{Fo}.\footnote{Foote studies the Prytz planimeter.} We identify all unit spheres $S^{n-1}$ along a curve $\Gamma(x)$ by parallel translations.

\begin{theorem} \label{hyper}
For all $x$, one has: $M_x\in O(n,1)$. 
\end{theorem}

\paragraph{\bf Proof.} Note first that the rear wheel's velocity is proj$_{r}v = (r\cdot v) r$, where $ r=AB $. Since $ M_x $ is the map of the sphere centered at the front wheel, we consider the moving  frame with the origin at the front wheel. This frame undergoes parallel translation as the wheel moves with its speed $v$. In the moving frame, the rear wheel's velocity is 
$ \omega (v) = (- v + (r\cdot v) r) \perp r $. We thus have a vector field on the sphere, and our map $ M_x $ is the time $x$--map of this vector field. It suffices therefore to show that this 
vector field corresponds to an element of the Lie algebra of $ o(n,1) $. 

The Lie algebra $o(n,1)$ consists of the matrices 
$$
C(M,v)=\begin{pmatrix}
M & v\\
v^* & 0
\end{pmatrix}
$$
where $M\in o(n)$ is an $n\times n$ skew-symmetric matrix and $v$ is an $n$-dimensional vector, 
and includes matrices of special form $ C(0,v)\:= C(v) $. We will show that these special matrices generate the vector field $ \omega (v) $  mentioned above. 
(As a side remark, the Lie algebra $o(n,1)$ is generated by its $n$-dimensional subspace $C(0, {\mathbb R} ^n)$). 

Let us compute the action of $C(v)$ on the unit sphere $S^{n-1}$. For a unit $n$-dimensional vector $r$, consider the point $(r,1)$ of the null cone at height $1$. Then
$$
\left( E + \eps C(v)\right) \begin{pmatrix} r\\1\end{pmatrix}= \begin{pmatrix}r+\eps v\\1+\eps r\cdot v\end{pmatrix} = 
k\begin{pmatrix} r-\eps \omega (v) \\1\end{pmatrix} + O( \eps^2),
$$
where $ k = (1 +\eps r\cdot v )$. 
 Thus $C(v)$ corresponds to the vector field $\omega (v)$ on the sphere,  and the result follows.
\proofend

\begin{remark} \label{constcurv}
{\rm It is quite likely that an analog of Theorem \ref{hyper} holds if $\R^n$ is replaced by either spherical or hyperbolic space. We do not dwell on it here. }
\end{remark}

\begin{remark} \label{snake}
{\rm It is interesting to point out possible connection with the so-called snake charmer algorithm \cite{H-R} in which the monodromy also takes values in the M\"obius group.}
\end{remark}

Now we consider corollaries of Theorem \ref{hyper} in the case $n=2$. Recall the classification of orientation preserving isometries of the hyperbolic plane: an elliptic isometry is a rotation about a point of $H^2$, and the corresponding map of the circle at infinity is conjugated to a rotation; a hyperbolic isometry has two fixed points at infinity, one exponentially attracting and another repelling; a parabolic isometry has a unique fixed point at infinity with derivative 1, see, e.g., \cite{Bea}. 

Let $\Gamma$, the trajectory of the front wheel, be closed. Then the monodromy map $M$ along $\Gamma$ is well-defined, up to conjugation; in particular, its type (elliptic, parabolic, hyperbolic) does not depend on the starting point.

The first corollary concerns the case when $M$ is hyperbolic. 

\begin{corollary} \label{closed} Let $M$ be hyperbolic, and 
let the trajectory of the rear wheel, $\gamma$, be a generic closed wave front. Then the trajectory of the front wheel, $\Gamma$, is also closed, and there exists a unique other closed trajectory of the rear wheel $\gamma^*$ with the same   front wheel trajectory $\Gamma$. The correspondence $\gamma \leftrightarrow \gamma^*$ is an involution. For a fixed orientation of $\Gamma$, one of the curves, $\gamma$ and $\gamma^*$, is exponentially stable and another exponentially unstable. The unstable curve $ \gamma$ is the closed path of the bike ridden backwards. 
\end{corollary}
  
 \paragraph{\bf Proof.} Since $\gamma$ is closed, the monodromy $M$ has a fixed point, and since $\gamma$ is generic, $M$ is hyperbolic. Then $M$ has another fixed point, corresponding to the closed trajectory $\gamma^*$. One of these fixed points is exponentially stable and another unstable. 
\proofend

\begin{figure}[hbtp]
\centering
\includegraphics[width=5.5in]{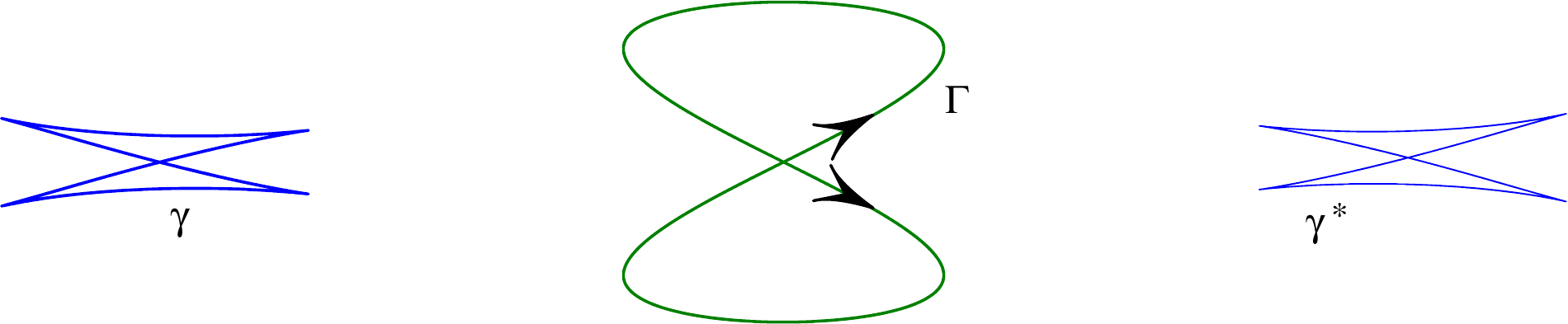}
    	\caption{The unstable curve is on the right. If the direction of traversal of figure eight is 	  reversed, the two curves exchange stability. }
    	\label{eight}
\end{figure}

\begin{figure}[hbtp]
\centering
\includegraphics[width=5.2in]{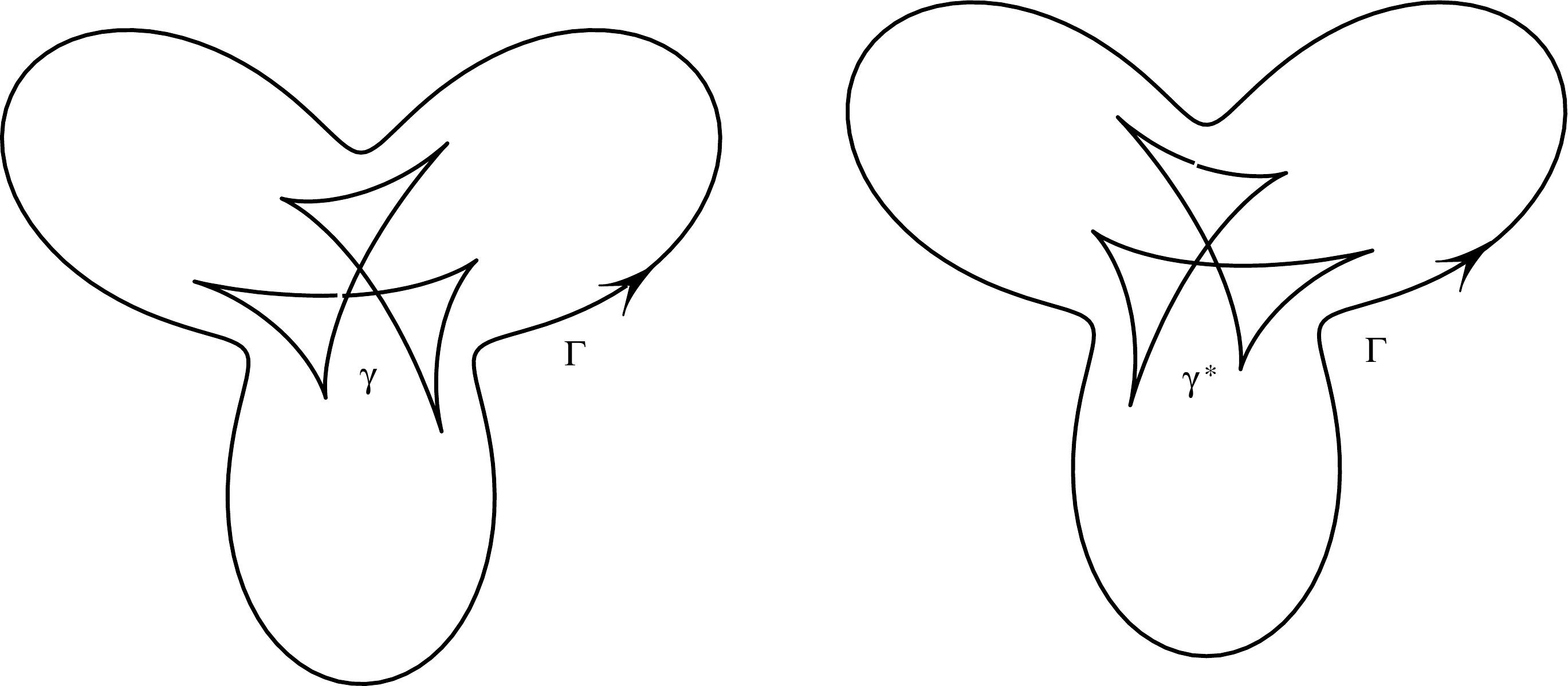}
    	\caption{The stable and the unstable rear trajectories for the sharmorck, just before the bifurcation when $\gamma$ and $ \gamma ^\ast $ coalesce. The shamrock is given by 
   $ x = r(t) \cos t, \ y =  r(t) \cos t $ with $ r = 0.94(1-0.5 \sin 3t).  $  }
    	\label{shmrck}
\end{figure}

Corollary \ref{closed} is illustrated by figures \ref{eight} and \ref{shmrck}.

 We precede  the next observation with a remark: for any  M\"obius map with two fixed points,   the derivatives at the two fixed points are reciprocal to each other. This, according to the next theorem, implies that $ \gamma $ and $ \gamma ^\ast $ have 
 the same length (up to a sign). 
 
 Let $\gamma$ be a closed wave front (the rear wheel track) and let $\Gamma=T(\gamma)$ be the front wheel track. Let $M$ be the monodromy of the curve $\Gamma$ and $L$ be the perimeter length of $\Gamma$.

 \begin{theorem} \label{hyperb}
Let  $M$ be hyperbolic or parabolic, and let $\gamma$ be the closed path of the rear wheel corresponding to a   fixed point $ \theta _0 $  of the M\"obius circle map $ \theta \mapsto M ( \theta ) $. Then 
\begin{equation} 
	 M ^\prime ( \theta _0 ) = e^{-\Le (\gamma )}.    
	 \label{eq:deriv}
\end{equation} 
 \end{theorem}
 
\begin{corollary}   If  $M$ is hyperbolic and $\gamma$ and $ \gamma ^\ast $ are the   rear tracks corresponding to the two fixed points then the curves $\gamma$ and $ \gamma ^\ast $ have equal lengths. 
\end{corollary}

\paragraph {Proof of the corollary.} For the fixed points $ \theta _0 $, $ \theta _0^\ast $ of any   M\"obius map one has  $ M ^\prime ( \theta _0 ) M ^\prime ( \theta _0^\ast ) =1 $, and the statement follows from Theorem \ref{hyperb}. \proofend

\begin{remark} \label{coi}
{\rm The case $\gamma=\gamma^*$ is quite interesting: this is  when one cannot tell which way the bicycle went from closed  tire tracks of the front and rear wheels, see Section \ref{intro}.
}
\end{remark}

\paragraph {\bf Proof of Theorem \ref{hyperb}.} Using the notation of Section \ref{prel}, consider equation (\ref{inverse}) (with $\ell=1$). This equation has an $L$-periodic solution $\alpha(x)$. Consider an infinitesimal perturbation     $ \alpha (x) + \eps \beta (x) $; the derivative of the monodromy map is given by 
$ M  ^\prime ( \theta _0) = \beta (L)/ \beta (0) $. But $\beta$ satisfies the linearized equation 
$ \beta ^\prime + \beta \cos \alpha = 0$, from which we find 
$$
	 M  ^\prime ( \theta _0) = \frac{\beta (L)}{\beta (0)} = e ^{- \int_{0}^{L} \cos \alpha (x) dx}.  
$$ 
But    $ \cos \alpha (x) $ is the speed of the rear wheel, and thus
$ \int_{0}^{L} \cos \alpha (x) dx  = \Le (\gamma )$. 
\proofend

\begin{remark} \label{id}
{\rm It is interesting that the monodromy may be identical, that is, there exist closed trajectories of the front wheel for which every trajectory of the rear wheel is closed. To construct such an example, let $\Gamma$ be  a small simple closed curve. Then the monodromy $M$ is elliptic, see analysis in \cite{Fo}. (This also follows from equation (\ref{inverse}): in the limit $\ell\to\infty$, the equation becomes $\alpha'(x)=\kappa(x)$, and since $\int \kappa(x)\ dx =2\pi$, the function $\alpha(x)$ cannot be periodic.) Slightly deforming $\Gamma$, if necessary, one may assume that $M$ is conjugated to a periodic rotation. Then, traversing $\Gamma$ an appropriate number of times, the monodromy becomes identical. }
\end{remark}

In contrast, if $ \Gamma$ is a closed immersed curve (not necessarily simple)
and $\ell$ is sufficiently small, one has a hyperbolic monodromy. Indeed, in the limit $\ell\to 0$, equation (\ref{inverse}) becomes $\sin\alpha=0$ and has two solutions $\alpha(x)=0$ and $\alpha(x)=\pi$, corresponding to the forward and backward tangent vectors to $\Gamma$. The two exponentially stable and unstable solutions survive for $\ell$ small enough.

As a limiting case of Theorem \ref{hyperb} for the parabolic monodromy we have the following.  
\begin{theorem} \label{parab}
$M$ is parabolic if and only if the total algebraic length of $\gamma$ is zero.
\end{theorem}

\paragraph{\bf Proof.} At the fixed point $ \theta _0 $ we have $ M ^\prime ( \theta _0 )  = 1$;  comparison with (\ref{eq:deriv}) shows  that $ \Le ( \gamma ) = 0 $. 
\proofend

\begin{corollary} \label{cusps}
In the parabolic case, the curve $\gamma$ has cusps.
\end{corollary}

\begin{figure}[hbtp]
\centering
\includegraphics[width=1.3in]{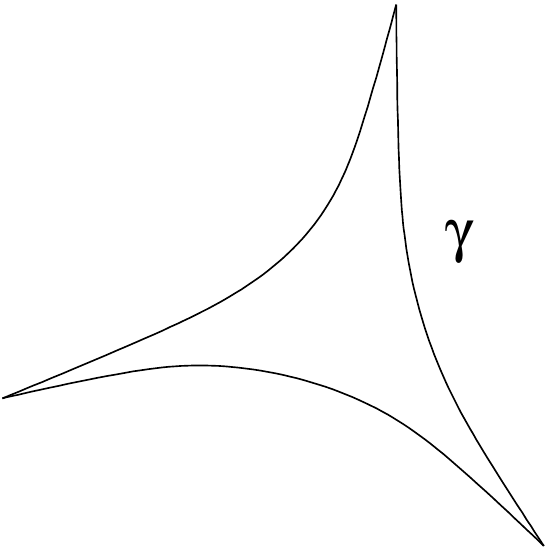}
\caption{Curve with turning number $\pi$}
\label{delta}
\end{figure}

An example of a wave front $\gamma$ yielding parabolic monodromy is depicted in figure \ref{delta}. The curve $\gamma$ has total turning number $\pi$, so for
$\Gamma$ to close up, one  traverses $\gamma$ twice. This ``doubled" front $\gamma$ obviously has zero total length.

\begin{figure}[thb]
   \centering
\includegraphics[width=5.2in]{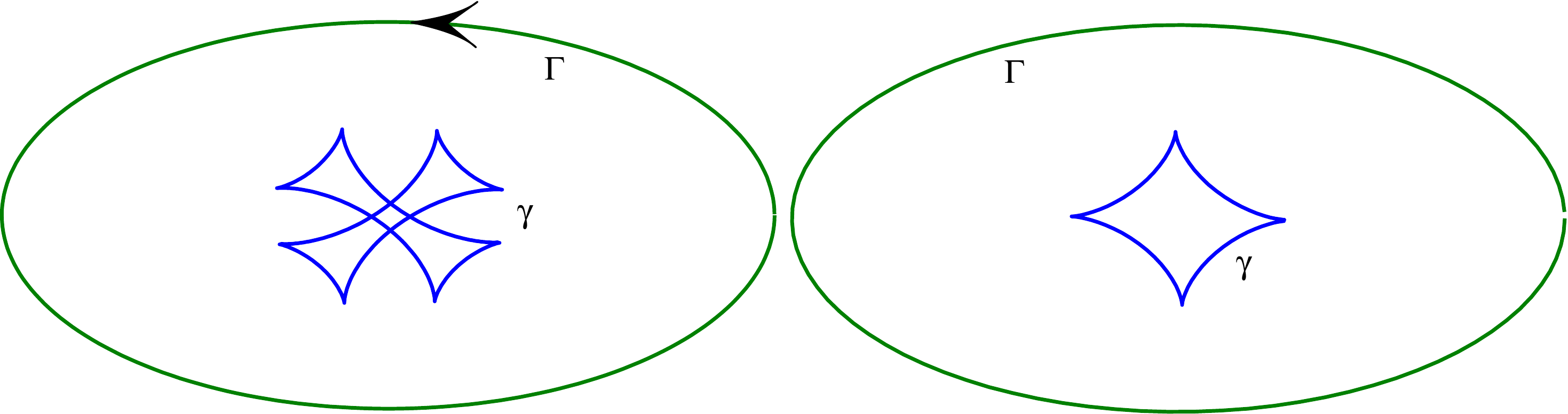}
\caption{A saddle--node bifurcation: the ellipse on  the left is (slightly) 
	larger.  The length of the coalesced curve on the right is zero  in accordance with Theorem \ref{parab}. In this particular case this is seen directly: the four arcs are congruent by symmetry, and their signs alternate. }
    	\label{fig:bifurcation}
\end{figure}

An example of the saddle--node bifurcation from the hyperbolic to the elliptic case, as the size of $ \Gamma $ decreases,  is shown in figure~\ref{fig:bifurcation}. 
 
\begin{remark} \label{iter}
{\rm Computation of the monodromy amounts to multiplying infinitely many $2\times 2$ matrices corresponding to infinitesimal arcs of the curve $\Gamma$ (if $\Gamma$ is a polygon, one has a finite products of hyperbolic elements in $SL(2, {\mathbb R} )$). A similar problem concerning the group of isometries of the sphere $SO(3)$  is treated in \cite{L1,L2}; we plan to extend this work to the group of isometries of the hyperbolic plane.
}
\end{remark}

\section{Proof of the Menzin conjecture} \label{Menzin}

\begin{theorem} \label{Menz}
 If $\Gamma$ is a closed convex curve bounding area greater than $\pi$ then the respective monodromy is hyperbolic.
\end{theorem}

\paragraph{\bf Proof.}  By approximation, we may assume that  $\Gamma$ is an oval, that is, a smooth closed strictly convex curve. We need to prove that if the monodromy $M$ is elliptic or parabolic then $\A(\Gamma)\leq \pi$. As we already mentioned, if $\Gamma$ is large enough, the monodromy $M$ is hyperbolic. Hence, if $M$ is elliptic, we can make $\Gamma$ larger (say, by homothety) and render $M$ parabolic. Therefore it suffices to prove that if $M$ is parabolic then $\A(\Gamma)\leq \pi$.

\begin{figure}[hbtp]
\centering
\includegraphics[width=3in]{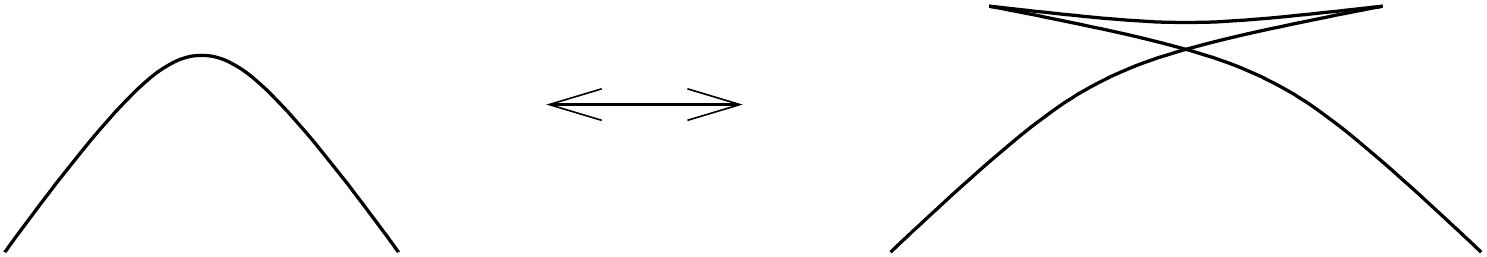}
\caption{Birth/death of a pair of cusps}
\label{tail}
\end{figure}

The proof is based on two observations: 
\begin{itemize} 
\item  $ \hbox {area}(\Gamma)=\hbox {area}(\gamma)+ \pi$, 
so that   $\A(\Gamma)\leq \pi$ is equivalent to $ \hbox{area} (\gamma )\leq 0 $,  and 
\item If $ \hbox {length}( \gamma ) =  0 $ then  $ \hbox{area} (\gamma )\leq 0 $.  
\end{itemize} 
We proceed with the detailed proof. 
Consider a generic one-parameter family of decreasing ovals $\Gamma_t,\ t\in [0,1]$, starting with a very large oval $\Gamma_0$ and ending with the given oval $\Gamma=\Gamma_1$ such that all the monodromies $M_t$ for $t\in[0,1)$ are hyperbolic. Therefore one has a family of wave fronts $\gamma_t$ (the closed trajectories of the rear wheel). Since $\Gamma_0$ is large enough, $\gamma_0$ is also an oval. The Legendrian liftings of the fronts $\gamma_t$ form a continuous family of immersed Legendrian curves in the space of contact elements. Therefore the Maslov index of $\gamma_1$ equals that of $\gamma_0$, that is, zero. Likewise, the rotation number $\rho(\gamma_1)$ equals one. The number of cusps may change in the family $\gamma_t$, see figure \ref{tail}.

The following holds due to the convexity of $\Gamma$. 
\begin{lemma} \label{infl}
The wave front has no inflections.
\end{lemma}

\begin{figure}[hbtp]
\centering
\includegraphics[width=3.5in]{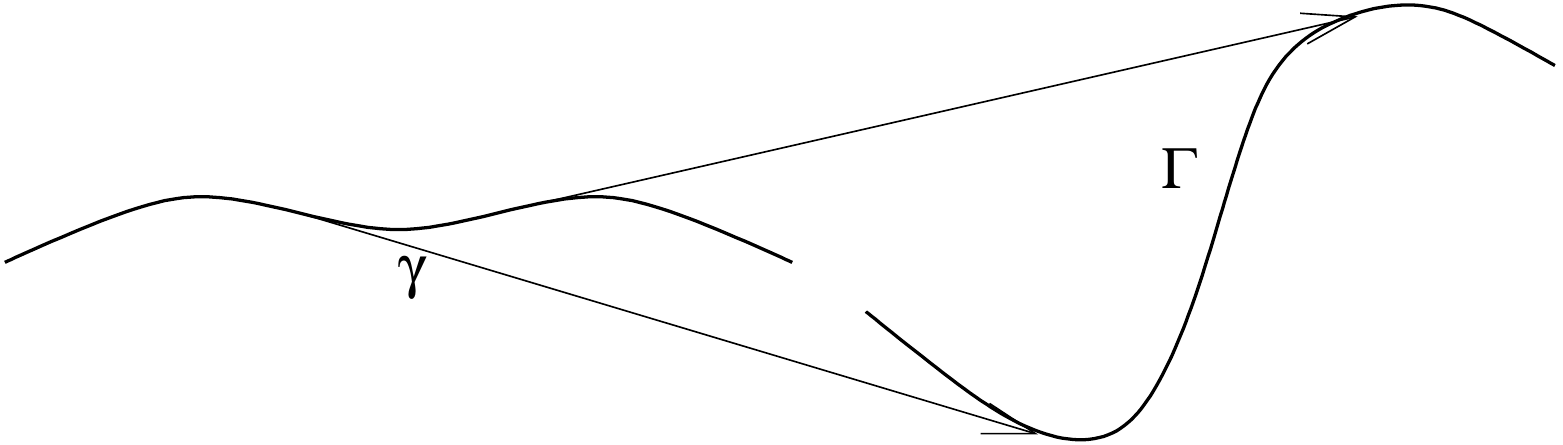}
\caption{Inflections of $\gamma$}
\label{dimple}
\end{figure}

\paragraph{\bf Proof.} Assume $\gamma_1$ has an inflection point. Note that $\gamma_0$ is convex. Let $\tau$ be the first value of parameter $t$ for which the curvature of $\gamma_t$ vanishes. Then, for $t$ slightly greater than $\tau$, the curve $\gamma_t$ has a ``dimple" and $\Gamma_t$ is not convex, see figure \ref{dimple}.
\proofend

Thus $\gamma_1$ is a wave front made of an even number of convex smooth arcs; the adjacent arcs form cusps. The total turning of the tangent direction to $\gamma_1$ is $2\pi$. The arcs are marked by $\pm$, the sign changes at each cusp. By Theorem \ref{parab}, the algebraic length of $\gamma_1$ vanishes: $\Le(\gamma_1)=0$.

Consider a smooth arc of $\gamma_1$ in the arclength parameterization; abusing notation, call this arc $\gamma_1(x)$. The respective arc of $\Gamma_1$ is
$\Gamma_1(x)=\gamma_1(x)+\sigma \gamma_1'(x)$
where $\sigma=\pm$ is the sign of the arc $\gamma_1$. Therefore
$\Gamma_1' = \gamma_1'+\sigma \gamma_1''$,
hence 
$$
\Gamma_1 \times \Gamma_1' = \gamma_1 \times \gamma_1' + \sigma \gamma_1 \times \gamma_1'' + \sigma^2 \gamma_1' \times \gamma_1''.
$$
Note that $\gamma_1 \times \gamma_1''=(\gamma_1 \times \gamma_1')'$ and that $ \gamma_1' \times \gamma_1''=k$, the curvature of $\gamma_1$.
Now we can find the area bounded by $\Gamma_1$:
\begin{equation} \label{area}
2\A(\Gamma_1)=\sum_i \left(\int \gamma_1(x) \times \gamma_1'(x)\ dx + \sigma_i \Delta_i (\gamma_1 \times \gamma_1') + \theta_i\right)
\end{equation}
where the sum is taken over the smooth arcs of $\gamma_1$, where $\sigma_i$ is the sign of $i$th arc, $\Delta_i$ is the difference of the momenta $\gamma_1 \times \gamma_1'$ at the end points of $i$th arc, and $\theta_i$ is the angle of turning of $i$th arc.  

Note that the sum of integrals in (\ref{area}) is $2\A(\gamma_1)$. Note also that $\sum \theta_i =2\pi$. Finally note that the terms $\Delta$ cancel out: $\sum \sigma_i \Delta_i (\gamma_1 \times \gamma_1') =0$. Therefore the inequality $\A(\Gamma_1)\leq \pi$ is equivalent to $\A(\gamma_1)\leq 0$.

To prove the latter inequality, let $p(\varphi)$ be the support function of the front $\gamma_1$ (the signed distance from the origin to the tangent line to $\gamma_1$ as a function of the direction of this line; see, e.g., \cite{San} for the theory of support functions). 
The support function exists because $\gamma_1$ is free from inflections and makes one full turn. One has the following formulas:
$$
\Le(\gamma_1)=\int_0^{2\pi} p(\varphi)\ d\varphi,\quad \A(\gamma_1)=\frac{1}{2} \int_0^{2\pi} (p^2(\varphi) - p'^2(\varphi) \ d\varphi.
$$
Thus we need to show that if 
$$
\int_0^{2\pi} p(\varphi)\ d\varphi =0\quad {\rm then}\quad \int_0^{2\pi} p^2(\varphi)\ d\varphi \leq \int_0^{2\pi} p'^2(\varphi)\ d\varphi.
$$
 But this is a well known Wirtinger inequality, which concludes the proof.
\proofend

\begin{remark} \label{iso}
{\rm Wirtinger inequality is intimately related with the isoperimetric inequality. Consider an oval $\gamma$ with area $A$ and perimeter length $L$. Consider the one-parameter family of equidistant fronts $\gamma_t$ inside the oval (that is, consider $\gamma$ as a source of light  propagating  inwards). The support function of $\gamma_t$ is that of $\gamma$ minus $t$. One has: 
$$\Le(\gamma_t)=L-2\pi t, \quad \A(\gamma_t)=A-Lt+\pi t^2.$$
 By the Wirtinger inequality, when $\Le(\gamma_t)=0$, one has $\A(\gamma_t)\leq 0$. Therefore if $t=L/2\pi$ then $A-Lt+\pi t^2\leq 0$, that is, $A\leq L^2/4\pi$, which is the isoperimetric inequality.
}
\end{remark}

\section{Oscillation of unicycle tracks} \label{osc}

Recall Finn's construction described in Section \ref{intro}.
Let $\gamma(t),\ t\in[0,L]$ be an arc length parameterized smooth curve in $\R^2$ such that the $\infty$-jets of $\gamma(t)$ coincide, for $t=0$ and $t=L$, with the $\infty$-jets of the $x$-axis at points $(0,0)$ and $(1,0)$, respectively. We use $\gamma$ as a ``seed" trajectory of the rear wheel of a bicycle. Then $\Gamma=T(\gamma)=\gamma+\gamma'$  is also tangent to the horizontal axis with all derivatives at its end points $(1,0)$ and $(2,0)$. Iterating this procedure yields a smooth infinite forward bicycle trajectory ${\cal T}$ such that the tracks of the rear and the front wheels coincide.
We shall study oscillation properties of ${\cal T}$. For starters, we note that the length of each new arc of ${\cal T}$ increases compared to the previous one.
 
 \begin{lemma} \label{incr}
The length of $\Gamma$ equals
$$
\int_0^L \sqrt{1+k^2(t)}\ dt >L
$$
where $k(t)=|\gamma''(t)|$ is the curvature  of $\gamma$. 
\end{lemma}

\paragraph{\bf Proof.} One has:
$$
\Gamma'(t)=\gamma'(t)+\gamma''(t),\quad |\Gamma'(t)|^2=1+|\gamma''(t)|^2,
$$
therefore the length of $\Gamma$ is
$$
\int_0^L |\Gamma'(t)|\ dt = \int_0^L \sqrt{1+k^2(t)}\ dt.
$$
\proofend

Denote by $Z(\gamma)$ the number of intersection points of the curve $\gamma(t),\ t\in (0,L)$ with the $x$-axis (we exclude the end points);  assume that $Z(\gamma)$ is finite.

\begin{proposition} \label{inter}
One has: $Z(\Gamma)> Z(\gamma)$. 
\end{proposition}

\begin{figure}[hbtp]
\centering
\includegraphics[width=3in]{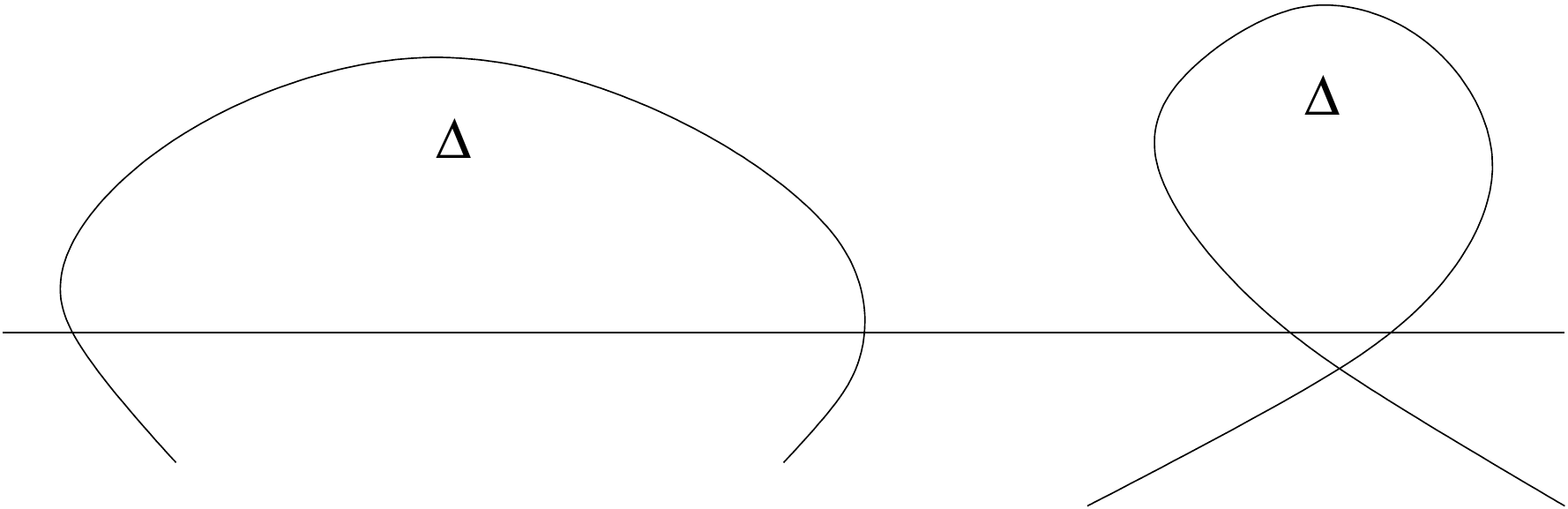}
\caption{Rolle theorem for curves}
\label{Rolle}
\end{figure}

\paragraph{\bf Proof.} Note that 
\begin{equation} \label{conj}
e^{-t} \left(e^t \gamma(t)\right)'=\Gamma(t).
\end{equation}
Let $Z(\gamma)=n$ and let $t_0=0<t_1<\dots<t_n<t_{n+1}=L$ be the consecutive moments of intersection  of $\gamma(t)$ with the $x$-axis. Then $t_i$ are also the consecutive moments of intersection  of the curve $\Delta(t):=e^t \gamma(t)$ with the $x$-axis. By a version of the Rolle theorem, see figure \ref{Rolle}, for each $i=0,1,\dots,n$, there is $t\in (t_i,t_{i+1})$ for which the curve $\Delta(t)$ has a horizontal tangent, i.e., the vector $\Delta'(t)$ is horizontal. It follows from (\ref{conj}) that $\Gamma(t)$ lies on the $x$-axis, and we are done.
\proofend

Consider the problem of extending the curve ${\cal T}$ backwards, that is, inverting the operator $T$. It turns out that usually $T$ can be inverted only finitely many times. Namely, one has the following corollary of Proposition \ref{inter}.

 \begin{corollary} \label{obstr}
 Let $\Gamma$ be a curve whose end points are unit distance apart and which is tangent, to all orders, at the end points to the $x$-axis. Let $Z(\Gamma)=n$. Then for no curve $\gamma$ whose end points are unit distance apart and which is tangent, to all orders, at the end points to the $x$-axis, one has $T^{n+1}(\gamma)=\Gamma$.
 \end{corollary}

Here is another oscillation property of the the curve ${\cal T}$. Let $E(\gamma)$ be the (finite) number of locally highest and lowest points of the curve $\gamma$. As before, $\Gamma=T(\gamma)$.

\begin{proposition} \label{extr}
One has: $E(\Gamma)> E(\gamma)$. 
\end{proposition}

\begin{figure}[hbtp]
\centering
\includegraphics[width=4in]{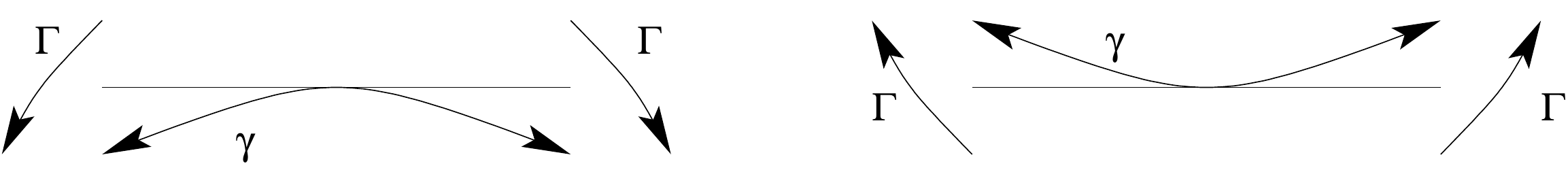}
\caption{Height extrema of the curve $\gamma$}
\label{extrem}
\end{figure}

\paragraph{\bf Proof.} At a locally highest point of $\gamma$, the curve $\Gamma$ has the downward direction, and at a locally lowest point -- the upward direction, see figure \ref{extrem}. It follows that, between consecutive locally highest and lowest points of $\gamma$, one has a locally lowest point of $\Gamma$, and between consecutive locally lowest and highest points of $\gamma$, one has a locally highest point of $\Gamma$. Considering the end points of $\gamma$ as local extrema of the height function as well, yields the result.
\proofend

\begin{conjecture} \label{conj1}
{\rm It follows from figure \ref{extrem} that the maximum height of $\Gamma$ is greater than that of $\gamma$, and likewise for the minimum height. We conjecture that the amplitude of the curve ${\cal T}$ is unbounded; in other words, unless $\gamma$ is a segment, ${\cal T}$ is not contained in any horizontal strip. We also conjecture that, unless $\gamma$ is a segment, ${\cal T}$ is not the graph of a function (i.e., one of the curves $T^n(\gamma)$ has a vertical tangent line) and, further, fails to be an embedded curve. One more conjecture: unless ${\cal T}$ is the horizontal axis, the curvature of ${\cal T}$ is unbounded.
}
\end{conjecture}

\paragraph{Configuration space of equilateral forward infinite linkages.}
 The construction of bicycle motion generating a single track can be interpreted as follows. Let ${\cal M}$ be the space of semi-infinite equilateral linkages $\{X=(x_0,x_1,x_2,\dots)\}$ where each $x_i$ is a point in the plane and $|x_i-x_{i+1}|=1$ for all $i$. Denote by $v_i$ the unit vector $x_i x_{i+1}$ and by $\alpha_{i}$ the angle between  $v_{i-1}$ and $v_{i}$. Let ${\cal M}_0$ be an open subset of ${\cal M}$ given by the condition $a_i \neq \pm \pi/2$ for all $i$. 
 
 \begin{figure}[hbtp]
\centering
\includegraphics[width=3.2in]{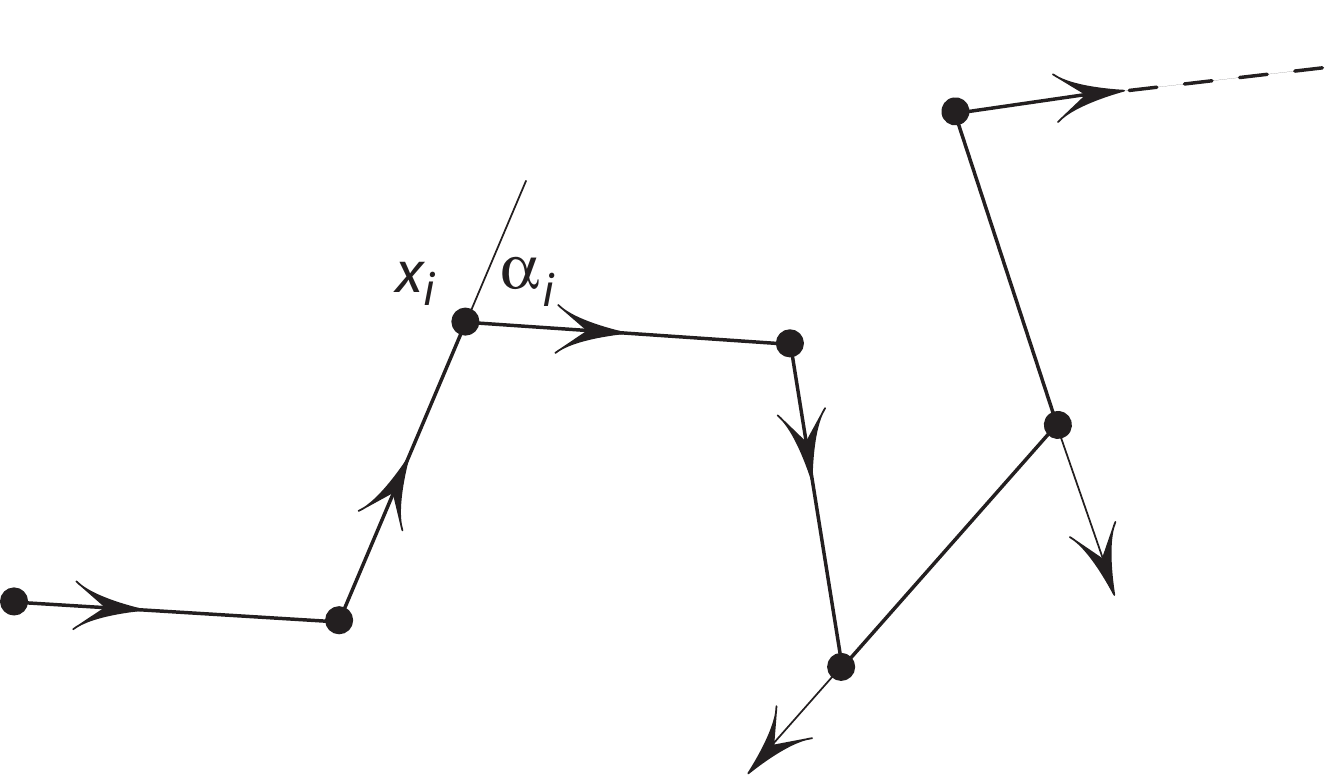}
    \caption{Note the change of direction when $ \alpha > \pi / 2  $. Only the direction of motion, and not the speeds, is indicated -- the latter becomes large for large values of $i$.  }
    \label{fig:linkage}
\end{figure} 

Consider the constraint on  ${\cal M}$ defined by the condition that the velocity of point $x_i$ is proportional to $v_i$. If $t_i$ is the speed of $x_i$ then the condition that  all links remain unit reads
\begin{equation} \label{stretch}
 t_i=t_{i+1} \cos\alpha_{i+1}
 \end{equation}
 for all $i$. On ${\cal M}_0$, where
$\cos\alpha_i\neq 0$, all the velocities are uniquely defined, up to a common factor, and one has a well defined field of directions $\xi$ which can be normalized to a vector field by setting $t_0=1$. If $\alpha_i=\pi/2$ for some $i$ then the speeds of all $x_j$ with $j<i$ must vanish; in particular, if $\alpha_i=\pi/2$ for infinitely many values of $i$ then such a configuration has no infinitesimal motions at all. See \cite{M-Z1,M-Z2} for  this non-holonomic system in relation to ``Monster Tower" and Goursat flags.

A forward bicycle motion generating a single track corresponds to  a solution to our system. The above described curve ${\cal T}$ yields an integral curve of the field $\xi$  in  ${\cal M}_0$: indeed,   $\alpha_i=\pi/2$ corresponds to a cusp of the trajectory of point $x_{i-1}$, whereas ${\cal T}$ is a smooth curve, as follows from  Lemma \ref{smooth}. 

The starting configuration $X$ of the Finn construction consists of non-negative integers on the horizontal axis, $x_i=(i,0)$, and one has a variety of integral curves of $\xi$ through $X\in {\cal M}_0$ (of which the simplest one is the uniform motion along the horizontal axis). Thus one has non-uniqueness of solutions of the differential equation describing the field $\xi$.

Finn's construction can be easily generalized as follows. Let $\delta$ be an infinite jet of a curve at point $x_0$. Consider the infinite jet  $T(\delta)$ at point $x_1=T(x_0)$, and let $\gamma$ be a curve smoothly interpolating between $\delta$ and $T(\delta)$. Then the concatenation of the curves $\gamma, T(\gamma), T^2(\gamma)$, etc., is a smooth unicycle track left by the bicycle motion with the seed curve $\gamma$.

The above construction provides a mapping $\Phi: J^{\infty} (x_0) \to {\cal M}_0$ from the space of infinite jets off curves at point $x_0$ to unit forward infinite linkages $\{(x_0,x_1,\dots)\}$.

\begin{proposition} \label{invert}
The mapping $\Phi$ is a bijection.
\end{proposition}

\paragraph{\bf Proof.} We construct the inverse map $\Psi: {\cal M}_0 \to J^{\infty} (x_0)$. Let $X=(x_0,x_1,\dots) \in {\cal M}_0$ and set: $C_i=\cos \alpha_i \neq 0$.
Then, according to (\ref{stretch}),
$$
t_0=1,\ \  t_k=\frac{1}{\Pi_{i=1}^k C_i},
$$
hence the speeds of all points are determined.  

We claim that, for each  $r\geq 1$, one has: $x_j^{(r)}=F_{j,r}(x_i,C_i)$ where $F$ is a polynomial in $x_i$ and a Laurent polynomial in $C_i$ for $i=0,1,\dots$.
This is proved by induction on $r$. For $r=1$, one has $x'_j=t_j (x_{j+1}-x_j)$. If $x_j^{(r)}=F_{j,r}(x_i,C_i)$ then
$$
x_j^{(r+1)}=\sum_i \frac{\partial F_{j,r}}{\partial x_i} x'_i + \frac{\partial F_{j,r}}{\partial C_i} C'_i.
$$
The induction step will be completed if we show that $C'_i$ is also a polynomial in $x_i$ and $C_i$. Indeed, $C_i = (x_i - x_{i-1})\cdot (x_{i+1}-x_i)$ and hence 
$$
C'_i = (t_i(x_{i+1}-x_i)-t_{i-1}(x_i - x_{i-1}))\cdot (x_{i+1}-x_i)  +
$$
$$
(x_i - x_{i-1})\cdot (t_{i+1}(x_{i+2}-x_{i+1}) -t_i(x_{i+1}-x_i)),
$$
as needed. 

In particular, $X$ determines all the derivatives $x_0^{(r)}$, that is, the infinite jet of a curve at $x_0$. This is $\Psi(X)$.
\proofend

We finish with another question:  is a straight line the only real analytic ``unicycle" trajectory?

\end{document}